\documentclass{amsart}
\usepackage{amsfonts}

\begin{document}
%\usepackage[inactive]{srcltx} % SRC Specials for DVI Searching

% Over-full v-boxes on even pages are due to the \v{c} in author's name
\vfuzz2pt % Don't report over-full v-boxes if over-edge is small

% THEOREM Environments ---------------------------------------------------
 \newtheorem{thm}{Theorem}[subsection]
 \newtheorem{cor}[thm]{Corollary}
 \newtheorem{lem}[subsection]{Lemma}
 \newtheorem{prop}[subsection]{Proposition}
 \theoremstyle{definition}
 \newtheorem{defn}[subsection]{Definition}
 \theoremstyle{remark}
 \newtheorem{rem}[subsection]{Remark}
 \numberwithin{equation}{section}
% MATH -------------------------------------------------------------------
\newcommand{\CC}{\mathbb{C}}
\newcommand{\KK}{\mathbb{K}}
\newcommand{\ZZ}{\mathbb{Z}}
\def\a{{\alpha}}

\def\b{{\beta}}

\def\d{{\delta}}

\def\g{{\gamma}}

\def\l{{\lambda}}

\def\gg{{\mathfrak g}}
\def\cal{\mathcal }

\title{infinitesimal deformations of the model $\ZZ_3$-filiform Lie algebra}

\author{R.M. Navarro}

\address{Rosa Mar{\'\i}a Navarro.\newline \indent
Dpto. de Matem{\'a}ticas, Universidad de Extremadura, C{\'a}ceres
(Spain)}

\email{rnavarro@unex.es}
%\thanks{This work has been partially supported by Junta de
%Extremadura-Consejer\'{\i}a de Infraestructuras y Desarrollo
%Tecnol\'ogico and by Feder (N. 3PR05A074)}

\begin{abstract}
In this work it is considered the vector space composed by the
infinitesimal deformations of the model $\ZZ_3$-filiform Lie
algebra $L^{n,m,p}$. By using these deformations all the
$\ZZ_3$-filiform Lie algebras can be obtained, hence the
importance of these deformations. The results obtained in this
work together to those obtained in \cite{dimension_color} and
\cite{erratum2}, leads to compute the total dimension of the
mentioned space of deformations.
\end{abstract}

%%% ----------------------------------------------------------------------
\maketitle
%%% ----------------------------------------------------------------------

\section{Introduction}

The concept of filiform Lie algebras was firstly introduced in
\cite{Vergne} by Vergne. This type of nilpotent Lie algebra has
important properties; in particular, every filiform Lie algebra
can be obtained by a deformation of the model filiform algebra
$L_n$. In the same way as filiform Lie algebras, all filiform Lie
superalgebras can be obtained by infinitesimal deformations of the
model Lie superalgebra $L^{n,m}$ \cite{Bor97}, \cite{JGP2},
\cite{JGP4} and \cite{JGP5}.

\

Continuing with the work of Vergne we have generalized the concept
and the propierties of the filiform Lie algebras into the theory
of color Lie superalgebras. Thus,   {\it filiform $G$-color Lie
superalgebras} and the model {\it filiform $G$-color Lie
superalgebra} were obtained in \cite{Filiformcolor}.

\

In the present paper the focus of interest are  color Lie
superalgebras with a $\ZZ_3$-grading vector space, i.e. $G=\ZZ_3$,
due to its physical applications
\cite{2},\cite{kerner2},\cite{kerner3}, \cite{5}, \cite{6} and
\cite{7}. Due to the fact that the one admissible commutation
factor for $\ZZ_3$ is exactly $\beta(g,h)=1 \ \forall g,h$,
$\ZZ_3$-color Lie superalgebras are indeed $\ZZ_3$-color Lie
algebras or $\ZZ_3$-graded Lie algebras. Thus, we have studied the
infinitesimal deformations of the model $\ZZ_3$-color Lie
superalgebra, i.e. the model $\ZZ_3$-filiform Lie algebra
$L^{n,m,p}$. By means of these deformations all $\ZZ_3$-filiform
Lie algebras can be obtained, hence the importance of these
deformations.

\

In \cite{dimension_color} and \cite{erratum2}, the authors
decomposed the space of these infinitesimal deformations, noted by
$Z^2(L;L)$, into six subspaces of deformations:

$$\begin{array}{l}Z^2(L;L)\cap \mathrm{Hom}(L_0 \wedge L_0,L_0)
\oplus Z^2(L;L)\cap \mathrm{Hom}(L_0 \wedge L_1,L_1)\oplus \\
Z^2(L;L)\cap \mathrm{Hom}(L_0 \wedge L_2,L_2) \oplus Z^2(L;L)\cap
\mathrm{Hom}(L_1 \wedge L_1,L_2) \oplus
\\  Z^2(L;L)\cap \mathrm{Hom}(L_1 \wedge L_2,L_0) \oplus Z^2(L;L)\cap
\mathrm{Hom}(L_2 \wedge L_2,L_1) \\
= A \oplus B \oplus C \oplus D \oplus E \oplus F
\end{array}$$

\

In the present paper it is  given a method that will allow  to
determine the dimension of the subspaces $A$, $B$ and $C$, giving
explicitly the total dimension of all of them (Theorems $1$, $2$
and $3$). This result together to those obtained in
\cite{dimension_color} and \cite{erratum2}, leads to obtain the
total dimension of the infinitesimal deformations of the model
$\ZZ_3$-filiform Lie algebra $L^{n,m,p}$ (Main Theorem).

\

We do assume that the reader is familiar
with the standard theory of Lie algebras. All the vector spaces
that appear in this paper (and thus, all the algebras) are assumed
to be  ${\mathbb F}$-vector spaces (${\mathbb F}=\CC$ or ${\mathbb
R}$) with finite dimension.

\section{Preliminaries}
The vector space $V$ is said to be $\ZZ_n-$graded if it admits a
decomposition in direct sum, $V=V_0 \oplus V_1 \oplus \cdots
V_{n-1}$. An element $X$ of $V$ is called homogeneous of degree
$\g$ ($deg(X)=d(X)=\g$), $\g \in \ZZ_n$, if it is an element of
$V_{\g}$.

\

Let $V=V_0 \oplus V_1 \oplus \cdots V_{n-1}$ and $W=W_0\oplus W_1
\oplus \cdots W_{n-1}$ be two graded vector spaces. A linear
mapping $f: V \longrightarrow W$ is said to be homogeneous of
degree $\g$ ($deg(f)=d(f)=\g$), $\g \in \ZZ_n$, if $ f(V_{\a})
\subset W_{\a + \g (mod \ n)}$ for all $\a \in \ZZ_n$. The mapping
$f$ is called a homomorphism of the $\ZZ_n-$graded vector space
$V$ into the $\ZZ_n-$graded vector space $W$ if $f$ is homogeneous
of degree 0. Now it is evident how we define an isomorphism or an
automorphism of $\ZZ_n-$graded vector spaces.

\

A superalgebra $\gg$ is just a $\ZZ_2-$graded algebra $\gg=\gg_0
\oplus \gg_1$. That is, if we denote by $[\  , \  ]$ the bracket
product of $\gg$, we have $[\gg_{\alpha}, \gg_{\beta}]\subset
\gg_{\alpha+\beta (mod 2)}$ for all $\alpha, \beta \in \ZZ_2$.

\

\begin{defn} \rm \label{defA} \cite{Scheunert} Let $\gg=\gg_0\oplus\gg_1$
be a superalgebra whose multiplication is denoted by the bracket
product [ , ]. We call $\gg$ a Lie superalgebra if the
multiplication satisfies the following identities:

 1. $[X,Y]=-(-1)^{\alpha \cdot \beta}[Y,X],\qquad \forall X\in \gg_{\alpha},
  \forall Y \in \gg_{\beta}$.

 2. $(-1)^{\g \cdot \alpha}[X,[Y,Z]]+(-1)^{\alpha \cdot \beta}[Y,[Z,X]]+
 (-1)^{\beta \cdot \g}[Z,[X,Y]]=0$
 \newline \indent \qquad for all $X\in \gg_{\alpha}, Y \in
 \gg_{\beta}, Z
 \in \gg_{\g}$ with  $\a, \b, \g  \in \ZZ_2$.

 \noindent Identity 2 is called the graded Jacobi identity and it will
 be denoted by $J_g(X,Y,Z)$.
\end{defn}

We observe that if $\gg=\gg_0\oplus\gg_1$ is a Lie superalgebra,
we have that $\gg_0$ is a Lie algebra and $\gg_1$ has the
structure of a $\gg_0-$module.

\

Color Lie (super)algebras can be seen as a direct generalization
of Lie (super)algebras. Indeed, the latter are defined through
antisymmetric (commutator) or symmetric (anticommutator) products,
although for the former the product is neither symmetric nor
antisymmetric and is defined by means of a commutation factor.
This commutation factor is equal to $\pm$ 1 for (super)Lie
algebras and more general for arbitrary color Lie (super)algebras.
As happened for Lie superalgebras, the basic tool to define color
Lie (super)algebras is a grading determined by an abelian group.

\begin{defn}
Let $G$ be an abelian group . A commutation factor $\beta$ is a
map
\newline $\beta$ : $G \times G \longrightarrow {\mathbb F} \setminus \{ 0
\}$, (${\mathbb F}=\CC$ or ${\mathbb R}$), satisfying the
following constraints:

\begin{itemize} \item[$(1)$] $\beta(g, h)\beta(h, g) = 1$ for all
$g, h \in G$ \item[$(2)$] $\beta(g, h + k) = \beta(g,
h)\beta(g,k)$ for all $g, h, k \in G$ \item[$(3)$] $\beta(g + h,
k) = \beta(g, k)\beta(h, k)$ for all $g, h, k \in G$
\end{itemize}

The definition above implies, in particular, the following
relations: $$\beta(0, g) = \beta(g, 0) = 1, \quad \beta(g, h) =
\beta(-h,g), \quad  \beta(g, g) =\pm 1 \quad \forall  g, h \in G
$$ where $0$ denotes the identity element of $G$. In particular,
fixing $g$ one element of $G$, the induced mapping $\beta_{g}:G
\longrightarrow {\mathbb F} \setminus \{ 0 \}$ defines a
homomorphism of groups.
\end{defn}

\begin{defn} Let $G$ be an abelian group and $\beta$ a commutation
factor. The (complex or real) $G-$graded algebra
$$L=\bigoplus_{g \in G}L{g}$$
with bracket product $[\  , \  ]$, is called a $(G,\beta)$-color
Lie superalgebra if for any $X \in L_g, \ Y \in L_h$, and $Z\in L$
we have

\begin{itemize}\item[(1)] $[X,Y]=-\beta(g,h)[Y,X]$ (anticommutative identity)
\item[(2)] $[[X,Y],Z]=[X,[Y,Z]]-\beta(g,h)[Y,[X,Z]]$ (Jacobi
identity)
\end{itemize}
\end{defn}

\begin{cor} Let $L=\bigoplus_{g \in G}L{g}$ be a
 $(G,\beta)$-color
Lie superalgebra. Then we have

\begin{itemize}\item[(1)] $L_0$ is a (complex or real) Lie algebra where $0$ denotes the identity element of
$G$.

\item[(2)] For all $g\in G \setminus \{ 0 \}$, $L_{g}$ is a
representation of $L_0$. If $X \in L_0$ and $Y \in L_g$, then
$[X,Y]$ denotes the action of $X$ on $Y$.
\end{itemize}
\end{cor}

\

 \noindent {\bf Examples.} For the particular case $G=\{0\}$, $L=L_0$ reduces to
a Lie algebra. If $G=\ZZ_2=\{0,1\}$ and $\b(1,1)=-1$ we have {\it
ordinary Lie superalgebras}, i.e. a {\it Lie superalgebra} is a
$(\ZZ_2,\b)$-color Lie superalgebra where $\b(i,j)=(-1)^{ij}$ for
all $i,j \in \ZZ_2$.

\begin{defn}

A representation of a $(G,\b)$-color Lie superalgebra is a mapping
$ \rho :L \longrightarrow End(V)$, where $V = \bigoplus_{g\in
G}V_g$ is a graded vector space such that
$$[\rho(X),\rho(Y)]=\rho(X)\rho(Y)-\b (g,h)\rho(Y)\rho(X)
$$

for all $X \in L_g, \ Y \in L_h$
\end{defn}

We observe that for all $g,h \in G$ we have $\rho(L_g)V_h
\subseteq V_{g+h}$, which implies that any $V_g$ has the structure
of a $L_0$-module. In particular considering the adjoint
representation $ad_{L}$ we have that every $L_g$ has the structure
of a $L_0$-module.

\

Two $(G,\b)$-color Lie superalgebras $L$ and $M$ are called {\bf
isomorphic} if there is a linear isomorphism $\varphi : L
\longrightarrow M$ such that $\varphi (L_g)=M_g$ for any $g\in G$
and also $\varphi ([x,y])=[\varphi(x),\varphi(y)]$ for any $x,y
\in L$.

\

Let $L=\bigoplus_{g \in G}L{g}$ be a $(G,\b)$-color Lie
superalgebra. The {\bf descending central sequence} of  $L$ is
defined by
$${\cal C}^0(L)=L, \quad
 {\cal C}^{k+1}(L)=[{\cal C}^k(L),L]\quad \forall k\geq 0$$ If
${\cal C}^k(L)=\{ 0 \}$ for some $k$, the $(G,\b)$-color Lie
superalgebra is called {\bf nilpotent}. The smallest integer $k$
such as ${\cal C}^k(L)=\{ 0 \}$ is called the {\bf nilindex} of
$L$.

\

Also, we are going to define some new {\bf descending sequences of
ideals}, see \cite{Filiformcolor}. Let $L=\bigoplus_{g \in G}L{g}$
be a
 $(G,\beta)$-color
Lie superalgebra. Then, we define the new descending sequences of
ideals ${\cal C}^{k}(L_0)$ (where $0$ denotes the identity element
of $G$) and ${\cal C}^{k}(L_{g})$ with $g\in G \setminus \{0\}$,
as follows:

 $${\cal C}^0(L_{0})=L_{0}, \quad {\cal C}^{k+1}(L_0)=[L_0, {\cal
 C}^k(L_0)], \quad
 k\geq 0$$
 and
 $${\cal C}^0(L_g)=L_g, \quad {\cal C}^{k+1}(L_g)=[L_0, {\cal
 C}^k(L_g)], \quad
 k\geq 0, \ g \in G \setminus \{0\}$$

Using the descending sequences of ideals defined above we give an
invariant of color Lie superalgebras called {\bf color-nilindex}.
We are going to particularize this definition for $G=\ZZ_3$.

\begin{defn} \cite{dimension_color}
 If $L=L_0\oplus L_1 \oplus L_2$ is a nilpotent
 $(\ZZ_3,\beta)$-color
Lie superalgebra, then $L$ has {\bf color-nilindex}
$(p_0,p_1,p_2)$, if the following conditions hold: $$ ({\cal
C}^{p_0-1}(L_0))({\cal C}^{p_1-1}(L_1))({\cal C}^{p_2-1}(L_2))\neq
0$$ and  $${\cal C}^{p_0}(\L_0)={\cal C}^{p_1}(\L_1)={\cal
C}^{p_2}(\L_2)=0$$
\end{defn}

\begin{defn} \cite{Filiformcolor} Let $L=\bigoplus_{g \in G}L{g}$ be a
 $(G,\beta)$-color
Lie superalgebra. $L_g$ is called a {\bf $L_0$-filiform module} if
there exists a decreasing subsequence of vectorial subspaces in
its underlying vectorial space $V$, $V=V_m \supset \dots \supset
V_1 \supset
 V_0$, with dimensions
$m,m-1,\dots 0$, respectively, $m>0$, and such that
$[L_0,V_{i+1}]= V_{i}$.
\end{defn}

\begin{rem} The definition of filiform module is also valid  for G-graded Lie algebras.
\end{rem}

\begin{defn} \cite{Filiformcolor} Let  $L=\bigoplus_{g \in G}L{g}$ be a
 $(G,\beta)$-color
Lie superalgebra. Then $L$ is a {\bf filiform color Lie
superalgebra} if the following conditions hold:
\begin{itemize}
\item[(1)] $L_0$ is a filiform Lie algebra where $0$ denotes the
identity element of $G$.

\item[(2)] $L_g$ has structure of $L_0$-filiform module, for all
$g \in G \setminus \{ 0 \}$
\end{itemize}
\end{defn}

\begin{defn} Let  $L=\bigoplus_{g \in G}L{g}$ be a
 $G$-graded
Lie algebra. Then $L$ is a $G$-{\bf filiform  Lie
algebra} if the following conditions hold:
\begin{itemize}
\item[(1)] $L_0$ is a filiform Lie algebra where $0$ denotes the
identity element of $G$.

\item[(2)] $L_g$ has structure of $L_0$-filiform module, for all
$g \in G \setminus \{ 0 \}$
\end{itemize}
\end{defn}
%For $G=\ZZ_3$ we give another equivalent definition for filiform
%color Lie superalgebras using the invariant called color-nilindex.
%
%
%\begin{defn} \cite{dimension_color} Any $(\ZZ_3,\b)$-color Lie superalgebra $L=L_0\oplus L_1 \oplus L_2$
%is a {\bf filiform color Lie superalgebra} if its color-nilindex
%is $(dim L_0-1,dim L_1,dim L_2)$.
%\end{defn}

It is not difficult to see that for $G=\ZZ_3$, there is only one
possibility for the commutation factor $\beta$, i. e.
$$\beta(g,h)=1  \qquad \forall \ g, \ h \ \in \ZZ_3=\{0,1,2\} $$

From now on we will consider this commutation factor and we will
write ``$\ZZ_3$-color" instead of ``$(\ZZ_3,\beta)$-color". We
will note by ${\cal L}^{n,m,p}$ the variety of all $\ZZ_3$-color
Lie superalgebras $L=L_0\oplus L_1 \oplus L_2$ with
$dim(L_0)=n+1$, $dim(L_1)=m$ and $dim(L_2)=p$. ${\cal N}^{n,m,p}$
will be the variety of all nilpotent $\ZZ_3$-color Lie
superalgebras  and ${\cal F}^{n,m,p}$ is the subset of ${\cal
N}^{n,m,p}$ composed of all filiform color Lie superalgebras.

\begin{rem}If $G=\ZZ_3$ then $\beta(g,h)=1  \ \forall \ g, \ h$. Thus, {\bf $\ZZ_3$-color Lie superalgebras are effectively
$\ZZ_3$-graded Lie algebras and filiform $\ZZ_3$-color Lie superalgebras are $\ZZ_3$-filiform Lie algebras}.
\end{rem}

In the particular case of $G=\ZZ_3$ the theorem of adapted basis
rest as follows for $L=L_0 \oplus \L_1 \oplus L_2 \in {\cal
F}^{n,m,p}$:

$$\left\{\begin{array}{ll}
[X_0,X_i]=X_{i+1},& 1\leq i \leq n-1,\\[1mm] [X_0,X_{n}]=0,&\\[1mm]
 [X_0,Y_j]=Y_{j+1},& 1\leq
j
\leq m-1,\\[1mm][X_0,Y_m]=0&\\[1mm]
[X_0,Z_k]=Z_{k+1},& 1\leq k
\leq p-1,\\[1mm][X_0,Z_p]=0.&
\end{array}\right.$$
with $\{X_0, X_1, \dots,$ $X_{n}\}$ a basis of  $L_0$,
$\{Y_1,\dots,Y_m\}$ a basis of $L_1$ and $\{Z_1,\dots,Z_p\}$ a
basis of $L_2$. The model $\ZZ_3$-filiform Lie algebra,
$L^{n,m,p}$, is the simplest  $\ZZ_3$-filiform Lie
algebra and it is defined in an adapted basis $\{X_0, X_1,
\dots, X_{n},$ $ Y_1, \dots, Y_m,$ $Z_1, \dots,Z_p \}$ by the
following non-null bracket products

$$L^{n,m,p}=\left\{\begin{array}{ll}
[X_0,X_i]=X_{i+1}, & 1 \leq i \leq n-1
\\[1mm]
[X_0,Y_j]=Y_{j+1},& 1 \leq j \leq m-1
\\ [1mm][X_0,Z_k]=Z_{k+1}& 1 \leq k \leq p-1
\end{array}\right.
$$

\section{cocycles and infinitesimal deformations}

Recall that a \textit{module} $V=V_0 \oplus V_1 \oplus V_2$ of the
$\ZZ_3$-color Lie superalgebra $L$ is a bilinear map of degree
$0$, $L\times V \to V$ satisfying
\[
   \forall~X\in L_g,~Y\in L_h ~v\in V:~~
     X(Y v)-Y(Xv)=[X,Y]v
\]
color Lie superalgebra cohomology is defined in the following
well-known way (see e.g. \cite{cohomology}): in particular, the
superspace of $q$-{\it dimensional cocycles} of the $\ZZ_3$-color
Lie superalgebra $L=L_0 \oplus L_1 \oplus L_2$ with coefficients
in the $L$-module $V=V_0 \oplus V_1 \oplus V_2$ will be given by
$$
 C^q(L;V)=\displaystyle{\bigoplus_{q_0+q_1+q_2=q} {\rm Hom} \left(\wedge^{q_0}
L_0 \otimes \wedge^{q_1} L_1 \otimes \wedge^{q_2} L_2, V \right)}
$$

This space is graded by $C^q(L;V)=C^q_0(L;V)\oplus C^q_1(L;V)
\oplus C^q_2(L;V)$ with
$$
  C^q_p(L;V)=
  \bigoplus_{\scriptsize \begin{array}{c}q_0+q_1+q_2=q \\
                        q_1+2q_2+p\equiv r \ mod \ 3
                          \end{array}}
  {\rm Hom} \left(\wedge^{q_0}
L_0 \otimes \wedge^{q_1} L_1 \otimes \wedge^{q_2} L_2, V_r \right)
$$

The {\it coboundary operator} $\delta^q: C^q(L;V) \longrightarrow
C^{q+1}(L;V)$, with  $\delta^{q+1} \circ \delta ^q=0$ is defined
in general, with $L$ an arbitrary $(G,\beta)$-color Lie
superalgebra and $V$ an $L$-module, by the following formula for
$q\geq 1$

$$\begin{array}{l}
 (\delta^q g)\big(A_0,A_1,\dots,A_q \big)=\\
   \hfill{\displaystyle \sum_{r=0}^{q} (-1)^{r} \beta (\gamma + \alpha_0 + \dots +
   \alpha_{r-1},\alpha_r)A_r \cdot
   g\big(A_0,\dots, \hat{A_r},\dots,A_q\big) }\\
%\qquad \qquad \qquad
   + \displaystyle \sum_{r < s} (-1)^{s} \beta (\alpha_{r+1} + \dots +
   \alpha_{s-1},\alpha_s)
   g\big(A_0,\dots, A_{r-1},[A_r,A_s],A_{r+1},\dots,\hat{A_s}, \dots
   A_q\big),
\end{array}$$

\noindent where $g \in C^q(L;V)$ of degree $\gamma$, and
$A_0,A_1,\dots,A_q \in L$ are homogeneous with degrees
$\alpha_0,\alpha_1,\dots,\alpha_q$ respectively. The sign $\hat{}$
indicates that the element below must be omitted and empty sums
(like $\alpha_0 + \dots + \alpha_{r-1}$ for $r=0$ and
$\alpha_{r+1} + \dots + \alpha_{s-1}$ for $s=r+1$) are set equal
to zero. In particular, for $q=2$ we obtain

$$\begin{array}{lll}
 (\delta^2 g)\big(A_0,A_1,A_2 \big)&=&
    \beta (\gamma , \alpha_0) A_0 \cdot
   g\big(A_1,A_2\big)-\beta (\gamma + \alpha_0,\alpha_1) A_1 \cdot
   g\big(A_0,A_2\big)+\\
%\qquad \qquad \qquad
&&   \beta (\gamma +\alpha_0+ \alpha_1, \alpha_2) A_2 \cdot
   g\big(A_0,A_1\big)\\
   &&-g([A_0,A_1],A_2)+\beta(\alpha_1,\alpha_2)g([A_0,A_2],A_1)+g(A_0,[A_1,A_2]).
\end{array}$$

Let $Z^q(L;V)$ denote the kernel of $\delta^q$ and let $B^q(L;V)$
denote the image of $\delta^{q-1}$, then we have that $B^q(L;V)
\subset Z^q(L;V)$. The elements of $Z^q(L;V)$ are called $q$-{\it
cocycles}, the elements of $B^q(L;V)$ are the $q$-{\it
coboundaries}. Thus, we can constuct the so-called {\it cohomology
groups}

$$ H^q(L;V)=Z^q(L;V)\left/ B^q(L;V) \right.$$
   $$ H^q_p(L;V)=Z^q_p(L;V)\left/ B^q_p(L;V) \right., if \ G=\ZZ_3 \ then \ p=0,1,2$$

 Two elements of $Z^q(L;V)$ are said to be {\it
cohomologous} if their residue classes modulo $B^q(L;V)$ coincide,
i.e., if their difference lies in $B^q(L;V)$.

\

We will focus our study in the 2-cocycles $Z^2_0(L^{n,m,p};
L^{n,m,p})$ with $L^{n,m,p}$ the model filiform $\ZZ_3$-color Lie
superalgebra. Thus $G=\ZZ_3$ and the only admissible  commutation
factor is exactly $\beta(g,h)=1$. Under all these restrictions the
condition that have to verify $\psi \in C^2_0(L^{n,m,p};
L^{n,m,p})$ to be a $2$-cocycle rests

$$\begin{array}{lll}
 (\delta^2 \psi)\big(A_0,A_1,A_2 \big)&=&
    [A_0,
   \psi \big(A_1,A_2\big)]- [A_1,
   \psi\big(A_0,A_2\big)]+\\
%\qquad \qquad \qquad
&&    [A_2,
   \psi \big(A_0,A_1\big)]-\psi([A_0,A_1],A_2)+\\
   && \psi([A_0,A_2],A_1)+ \psi(A_0,[A_1,A_2])=0
\end{array}$$

\noindent for all $A_0, A_1, A_2  \in L^{n,m,p}$. We observe that
$L^{n,m,p}$ has the structure of a $L^{n,m,p}$-module via the
adjoint representation.

\

We consider an homogeneous basis of $L^{n,m,p}=L_0 \oplus L_1
\oplus L_2$, in particular an adapted basis $\{X_0, X_1, \dots,$
$X_{n}, Y_1,\dots,Y_m,Z_1,\dots,Z_p \}$ with $\{X_0, X_1, \dots,$
$X_{n}\}$ a basis of $L_0$, $\{Y_1,\dots,Y_m\}$ a basis of $L_1$
and $\{Z_1,\dots,Z_p\}$ a basis of $L_2$.

Under these conditions we have the following lemma.

\begin{lem}\label{10conditions} \cite{dimension_color}, \cite{erratum2} Let $\psi$ be such that $\psi \in C^2_0(L^{n,m,p};L^{n,m,p})$, then $\psi$ is a $2$-cocycle, $\psi \in
Z^2_0(L^{n,m,p};L^{n,m,p})$, iff the $10$ conditions below hold
for all $X_i, X_j, X_k \in L_0$, $Y_i, Y_j, Y_k \in L_1$ and $Z_i,
Z_j, Z_k \in L_2$

\begin{itemize}
\item[$(1)$]
$[X_i,\psi(X_j,X_k)]-[X_j,\psi(X_i,X_k)]+[X_k,\psi(X_i,X_j)]-\psi([X_i,X_j],X_k)+$\\
$\psi([X_i,X_k],X_j)+\psi(X_i,[X_j,X_k])=0$
 \item[$(2)$]
$[X_i,\psi(X_j,Y_k)]-[X_j,\psi(X_i,Y_k)]+[Y_k,\psi(X_i,X_j)]-\psi([X_i,X_j],Y_k)+$\\
$\psi([X_i,Y_k],X_j)+\psi(X_i,[X_j,Y_k])=0$
 \item[$(3)$]
$[X_i,\psi(X_j,Z_k)]-[X_j,\psi(X_i,Z_k)]+[Z_k,\psi(X_i,X_j)]-\psi([X_i,X_j],Z_k)+$\\
$\psi([X_i,Z_k],X_j)+\psi(X_i,[X_j,Z_k])=0$
 \item[$(4)$]
$[X_i,\psi(Y_j,Y_k)]-[Y_j,\psi(X_i,Y_k)]+[Y_k,\psi(X_i,Y_j)]-\psi([X_i,Y_j],Y_k)+$\\
$\psi([X_i,Y_k],Y_j)+\psi(X_i,[Y_j,Y_k])=0$
 \item[$(5)$]
$[X_i,\psi(Y_j,Z_k)]-[Y_j,\psi(X_i,Z_k)]+[Z_k,\psi(X_i,Y_j)]-\psi([X_i,Y_j],Z_k)+$\\
$\psi([X_i,Z_k],Y_j)+\psi(X_i,[Y_j,Z_k])=0$
 \item[$(6)$]
$[X_i,\psi(Z_j,Z_k)]-[Z_j,\psi(X_i,Z_k)]+[Z_k,\psi(X_i,Z_j)]-\psi([X_i,Z_j],Z_k)+$\\
$\psi([X_i,Z_k],Z_j)+\psi(X_i,[Z_j,Z_k])=0$
 \item[$(7)$]
$[Y_i,\psi(Y_j,Y_k)]-[Y_j,\psi(Y_i,Y_k)]+[Y_k,\psi(Y_i,Y_j)]-\psi([Y_i,Y_j],Y_k)+$\\
$\psi([Y_i,Y_k],Y_j)+\psi(Y_i,[Y_j,Y_k])=0$
 \item[$(8)$]
$[Y_i,\psi(Y_j,Z_k)]-[Y_j,\psi(Y_i,Z_k)]+[Z_k,\psi(Y_i,Y_j)]-\psi([Y_i,Y_j],Z_k)+$\\
$\psi([Y_i,Z_k],Y_j)+\psi(Y_i,[Y_j,Z_k])=0$
 \item[$(9)$]
$[Y_i,\psi(Z_j,Z_k)]-[Z_j,\psi(Y_i,Z_k)]+[Z_k,\psi(Y_i,Z_j)]-\psi([Y_i,Z_j],Z_k)+$\\
$\psi([Y_i,Z_k],Z_j)+\psi(Y_i,[Z_j,Z_k])=0$
 \item[$(10)$]
$[Z_i,\psi(Z_j,Z_k)]-[Z_j,\psi(Z_i,Z_k)]+[Z_k,\psi(Z_i,Z_j)]-\psi([Z_i,Z_j],Z_k)+$\\
$\psi([Z_i,Z_k],Z_j)+\psi(Z_i,[Z_j,Z_k])=0$
\end{itemize}
\end{lem}

\begin{prop} \cite{dimension_color} $\psi$ is an infinitesimal deformation of $L^{n,m,p}$
iff $\psi$ is a $2$-cocycle of degree $0$, $\psi \in
Z^2_0(L^{n,m,p};L^{n,m,p})$.
\end{prop}

\begin{thm}\cite{Filiformcolor} {\it (1) Any filiform $(G,\beta)$-color Lie superalgebra law
$\mu$ is isomorphic to $\mu_0 + \varphi$ where $\mu_0$ is the law
of the model filiform $(G,\beta)$-color Lie superalgebra and
$\varphi$ is an infinitesimal deformation of $\mu_0$ verifying
that $\varphi(X_0,X)=0$ for all $X \in L$, with $X_0$ the
characteristic vector of model one.}

{\it (2) Conversely, if $\varphi$ is an infinitesimal deformation
of a model filiform $(G,\beta)$-color Lie superalgebra law $\mu_0$
with $\varphi(X_0,X)=0$ for all $X \in L$, then the law $\mu_0 +
\varphi$ is a filiform $(G,\beta)$-color Lie superalgebra law iff
$\varphi \circ \varphi =0$.}
\end{thm}

Thus, any  $\ZZ_3$-filiform Lie algebra (filiform $\ZZ_3$-color Lie superalgebra) will  be   a
linear deformation of the model $\ZZ_3$-filiform Lie algebra (the model $\ZZ_3$-color Lie
superalgebra), i.e. $L^{n,m,p}$ is the model $\ZZ_3$-filiform
Lie algebra an another arbitrary $\ZZ_3$-filiform Lie
algebra will be equal to $L^{n,m,p}+\varphi$, with $\varphi$
an infinitesimal deformation of $L^{n,m,p}$. Hence the importance
of these deformations. So, in order to determine all the
$\ZZ_3$-filiform Lie algebras it is only necessary to compute
the infinitesimal deformations or so called $2$-cocycles of degree
$0$, that vanish on the characteristic vector $X_0$. Thanks to the
following lemma these infinitesimal deformations will can be
decomposed into $6$ subspaces.

\begin{lem} \cite{dimension_color}, \cite{erratum2} Let $Z^2(L;L)$ be the $2$-cocycles $Z^2_0(L^{n,m,p};L^{n,m,p})$
that vanish on the characteristic vector $X_0$. Then, $Z^2(L;L)$
can be divided into six subspaces, i.e. if $L^{n,m,p}=L=L_0\oplus
L_1 \oplus L_2$ we will have that

\

$\begin{array}{rl}Z^2(L;L)=&Z^2(L;L)\cap \mathrm{Hom}(L_0 \wedge
L_0,L_0) \oplus Z^2(L;L)\cap \mathrm{Hom}(L_0 \wedge
L_1,L_1)\oplus \\& Z^2(L;L)\cap \mathrm{Hom}(L_0 \wedge L_2,L_2)
\oplus Z^2(L;L)\cap \mathrm{Hom}(L_1 \wedge L_1,L_2) \oplus
\\ & Z^2(L;L)\cap \mathrm{Hom}(L_1 \wedge L_2,L_0) \oplus Z^2(L;L)\cap
\mathrm{Hom}(L_2 \wedge L_2,L_1) \\
=& A \oplus B \oplus C \oplus D \oplus E \oplus F
\end{array}$

\end{lem}

In order to obtain the dimension of $A$, $B$ and $C$ we are going
to adapt {\it the $\mathfrak{sl}_2(\CC)$-module method} that we
have already used for Lie superalgebras \cite{Bor97}, \cite{JGP2},
\cite{JGP4} and for color Lie superalgebras
\cite{dimension_color}, \cite{erratum2}. Next, we will do it
explicitly for $A=Z^2(L;L)\cap \mathrm{Hom}(L_0 \wedge L_0,L_0) $.

\section{Dimension of $A=Z^2(L;L)\cap \mathrm{Hom}(L_0 \wedge
L_0,L_0) $}

\

In general, any cocycle $a \in Z^2(L;L)\cap \mathrm{Hom}(L_0
\wedge L_0,L_0)$ will be any skew-symmetric bilinear map from $L_0
\wedge L_0$ to $L_0$ such that:

\

\begin{itemize}
\item[$(1)$]
$[X_i,a(X_j,X_k)]-[X_j,a(X_i,X_k)]+[X_k,a(X_i,X_j)]-a([X_i,X_j],X_k)+$\\
$a([X_i,X_k],X_j)+a(X_i,[X_j,X_k])=0$ \hfill{ $\forall \
X_i,X_j,X_k \in L_0$}
\end{itemize}

\

\noindent with $a(X_0,X)=0$ $\forall X \in L$. As $X_0 \notin$ Im
$a$ and taking into account the bracket products of $L$ then the
equation $(1)$ can be rewritten as follows

\begin{equation} \label{cocycle}
[X_0,a(X_j,X_k)]-a([X_0,X_j],X_k)-a (X_j,[X_0,X_k])=0,\ 1 \leq j <
k \leq n
\end{equation}

\

In order to obtain the dimension of the space of cocycles for $A$
we apply an adaptation of the $\mathfrak{sl}(2,\CC)$-module Method
that we used in \cite{dimension_color}.

\

Recall the following well-known facts about the Lie algebra
$\mathfrak{sl}(2,\CC)$ and its finite-dimensional modules, see
e.g. \cite{Bourbaki7}, \cite{Humphreys}:

$\mathfrak{ sl}(2,\CC)=<X_{-},H,X_{+}>$ with the following
commutation relations:
$$
   \left\{\begin{array}{l}
      [X_{+},X_{-}]=H \\[1mm] [H,X_{+}]=2X_{+}, \\[1mm]
          [H,X_{-}]=-2X_{-}.
          \end{array}\right.
$$
Let $V$ be a $n$-dimensional  $\mathfrak{ sl}(2,\CC)$-module,
$V=<e_1,\dots,e_n>$. Then, up to isomorphism there exists a unique
structure of an
 irreducible $\mathfrak{ sl}(2,\CC)$-module in $V$ given in a basis
$\{e_1,\ldots,e_n\}$ as follows \cite{Bourbaki7}:
$$
   \left\{ \begin{array}{ll}
     X_{+}\cdot e_i=e_{i+1}, & 1 \leq i \leq n-1,\\[1mm]
            X_{+}\cdot e_n=0, &                  \\[1mm]
   H \cdot e_i=(-n+2i-1)e_i,  & 1 \leq i \leq n.
           \end{array}\right.
$$

It is easy to see that $e_n$ is the maximal vector of $V$ and its
weight, called the highest weight of $V$, is equal to $n-1$.

Let $W_0,W_1,\dots,W_k$ be $\mathfrak{ sl}(2,\CC)$-modules, then
the space $\mathrm{Hom}(\otimes_{i=1}^k W_i,W_0)$ is a
$\mathfrak{sl}(2,\CC)$-module in the following natural manner:
$$
  (\xi \cdot \varphi)(x_1,\dots,x_k)=\xi \cdot \varphi(x_1,\dots,x_k)-
   \sum_{i=1}^{k}
    \varphi (x_1,\dots,\xi \cdot x_i,x_{i+1},\dots,x_n)
$$
with $\xi \in \mathfrak{ sl}(2,\CC)$ and $\varphi
\in\mathrm{Hom}(\otimes_{i=1}^k W_i,W_0)$. In particular, if $k=2$
and $W_0=W_1=W_2=V_0$, then
$$
    (\xi \cdot \varphi)(x_1,x_2)=\xi \cdot \varphi(x_1,x_2)-
 \varphi (\xi \cdot x_1,x_2)- \varphi (x_1,\xi \cdot x_2).
$$

An element $\varphi \in\mathrm{Hom}(V_0\otimes V_0,V_0)$ is said
to be invariant if
%$$
   $X_{+}\cdot \varphi =0$,
%$$
that is
\begin{equation} \label{maximal}
      X_{+} \cdot \varphi(x_1,x_2)-
        \varphi (X_{+} \cdot x_1,x_2)- \varphi (x_1,X_{+} \cdot x_2)=0,
         \quad \forall x_1,x_2 \in V.
\end{equation}

Note that $\varphi \in\mathrm{Hom}( V_0\otimes V_0,V_0)$ is
invariant if and only if $\varphi$ is a maximal vector.

\

 We are going to
consider the structure of irreducible
$\mathfrak{sl}(2,\CC)$-module in $V_0=<X_1,\dots,X_n>=L_{0}/\CC
X_0$, thus in particular:
$$
   \left\{ \begin{array}{ll}
      X_{+}\cdot X_i=X_{i+1}, & 1 \leq i \leq n-1,\\[1mm]
              X_{+}\cdot X_n=0,&
             \end{array}\right.
$$

Next, we identify the multiplication of $X_{+}$ and $X_i$ in the
$\mathfrak{ sl}(2,\CC)$-module  $V_0=<X_1,\dots,X_n>$, with the
bracket $[X_0,X_i]$ in $L_{0}$ and thanks to these
identifications, the expressions (\ref{cocycle}) and
(\ref{maximal}) are equivalent. Thus  we have the following
result:

\begin{prop} Any skew-symmetric bilinear map $\varphi$,
$\varphi: V_0 \wedge V_0 \longrightarrow V_0$ will be an element
of the space of cocycles $A$ if and only if $\varphi$ is a maximal
vector of the $\mathfrak{ sl}(2,\CC)$-module
$\mathrm{Hom}(V_0\wedge V_0,V_0)$, with $V_0=\langle
X_1,\ldots,X_n\rangle$.
\end{prop}
%\begin{proof} If $\varphi \in C$, then $Im \varphi \subset
%V_0$. In fact, the ideal $V_0$ is equal to its own centralizer in $L^{n,m}_0$.
%Let $x\in V_0$ and $y,y'\in L^{n,m}_1$. Thanks to the cocycle
%identity we get
%\[
% [x,\varphi(y,y')]=\varphi(y,[x,y'])+\varphi(y',[x,y])=0
%  \quad \forall~x\in V_0; y,y' \in L^{n,m}_1.
%\]
%since $[V_0,V_1]=\{0\}$. Hence $\varphi(y,y')$ centralizes
%$V_0$, and is thus an element of $V_0$.
%%if $X_0 \in Im \varphi$ it can be supposed that
%%there exist $i,j$ such that $\varphi(Y_i,Y_j)=X_0$, but it
%%constitutes a contradiction with the cocycle equation
%%$$
%%  [x,\varphi(y,z)]-\varphi(z,[x,y])-\varphi(y,[x,z])=0,
%%\quad \forall x\in L^{n,m}_0; y,z \in L^{n,m}_1.
%%$$
%%for $x=X_1$, $y=Y_i$ and $z=Y_j$.
%The same cocyle equation for $x$ replaced by $X_0$, i.e.
%\[
%  [X_0,\varphi(y,y')]-\varphi(y,[X_0,y'])-\varphi(y',[X_0,y])=0,
%\]
%shows that $\varphi\in \mathrm{Hom}(S^2 V_1,V_0)$ is invariant
%which proves the proposition and the equation (\ref{cocycle}).
%\end{proof}

\begin{cor} As each irreducible $\mathfrak{sl}(2,\CC)$-module has (up to nonzero scalar
multiples) a unique maximal vector, then the dimension of the
space of cocycles $A$ is equal to the number of summands of any
decomposition of $\mathrm {Hom}(V_0\wedge V_0,V_0)$ into the
direct sum of irreducible $\mathfrak{sl}(2,\CC)$-modules.
\end{cor}

We use the fact that each irreducible module contains either a
unique (up to scalar multiples) vector of weight $0$ (in case the
dimension of the irreducible module is odd) or a unique (up to
scalar multiples) vector of weight $1$ (in case the dimension of
the irreducible module is even). We therefore have

%Thanks to the symmetric structure of the weights, instead of to
%sum the maximal vectors it is possible, and easier, to sum the
%vectors of weight 0 or 1.

\begin{cor}\label{cor2} The dimension of the
space of cocycles $A$ is equal to the dimension of the subspace of
{\rm Hom} $(V_0 \wedge V_0,V_0)$ spanned by the vectors of weight
0 or 1.
\end{cor}

At this point,  we are going to apply the
$\mathfrak{sl}(2,\CC)$-module method aforementioned in order to
obtain the dimension of the space of cocycles $A$.

\

We consider a natural basis $\mathcal{B}$ of $\mathrm{Hom}(V_0
\wedge V_0,V_0)$
%$$
%   {\mathcal B}=\{\varphi_{i,j}^s\}
%$$
%with
consisting of the following maps:
$$
   \varphi_{i,j}^s(X_k,X_l)=\left\{
         \begin{array}{ll}
           X_s & \mbox{ if }(i,j)=(k,l)\\
             0 & \mbox{ in all other cases}
         \end{array}
\right.
$$

\noindent where  $1\leq i,j,k,l,s \leq n$, with $i \neq j$ and
$\varphi_{i,j}^s=- \varphi_{j,i}^s$.

\

 Thanks to
Corollary \ref{cor2} it will be enough to find the basis vectors
$\varphi_{i,j}^s$ with weight $0$ or $1$. The weight of an element
$\varphi_{i,j}^s$ (with respect to $H$) is
$$
  \lambda(\varphi_{i,j}^s)
    =\lambda(X_s)-\lambda(X_i)-\lambda(X_j)=n+2(s-i-j)+1.
$$

In fact,
$$
    \begin{array}{ll}
(H \cdot \varphi_{i,j}^s)(X_i,X_j) & =H \cdot
\varphi_{i,j}^s(X_i,X_j)-\varphi_{i,j}^s(H \cdot
X_i,X_j)-\varphi_{i,j}^s(X_i,H \cdot X_j) \\ \\
& = H\cdot X_s
-\varphi_{i,j}^s((-n-1+2i)X_i,X_j)-\varphi_{i,j}^s(X_i,(-n-1+2j)X_j)\\
\\
&= (-n-1+2s)X_s-(-n-1+2i)X_s-(-n-1+2j)X_s \\ \\
& = [n+2(s-i-j)+1]X_s
    \end{array}
$$
%We are going to introduce a simpler weight,  it corresponds to the
%action of the diagonalizable derivation $d$, $d \in \mathrm{Der}
%L_0$,
% defined by:
%$$
%   d(X_0)=X_0,\ d(X_i)=iX_i; \quad 1\leq i \leq n.
%$$
%This weight will be denoted by $p(\varphi)$. We have that
%$$
%    p(\varphi_{i,j}^s)=s-i-j.
%$$
%We have the following relationships between the two weights:
%$$
%   \begin{array}{l}
%            \lambda(\varphi)=2p(\varphi)+n+1,\\
%                  p(\varphi)=\frac{1}{2}(\lambda(\varphi)-n-1).
%\end{array}
%$$

We observe that if $n$ is even then $\lambda(\varphi)$ is odd, and
if $n$ is odd then $\lambda(\varphi)$ is even. So, if $n$ is even
it will be sufficient to find the elements $\varphi_{i,j}^s$ with
weight $1$ and if $n$ is odd it will be sufficient to find those
of them with weight $0$.

We can consider the three sequences that correspond with the
weights of $V=<X_1,X_2,\dots,X_{n-1},X_n>$ in order to find the
elements with weight $0$ or $1$:
$$
  -n+1,-n+3,\dots,n-3,n-1;
$$
$$
  -n+1,-n+3,\dots,n-3,n-1;
$$
$$
  -n+1,-n+3,\dots,n-3,n-1.
$$

\noindent and we have to count the number of all possibilities to
obtain $1$ (if $n$ is even) or $0$ (if $n$ is odd). Remember that
$\lambda(\varphi_{i,j}^s)=\lambda(X_s)-\lambda(X_i)-\lambda(X_j)$,
where $\lambda(X_s)$ belongs to the last sequence, and
$\lambda(X_i)$, $\lambda(X_j)$ belong to the first and second
sequences respectively. For example, if $n$ is odd, we have to
obtain $0$, so we can fix an element (a weight) of the last
sequence and then count the possibilities to sum the same quantity
between the two first sequences. Taking into account the
skew-symmetry of $\varphi_{i,j}^s$, that is
$\varphi_{i,j}^s=-\varphi_{j,i}^s$ and $i \neq j$, and repeating
the above reasoning for all the elements of the last sequence we
obtain the following theorem:

\

\noindent {\bf Theorem 1.} Let $Z^2(L;L)$ be the $2$-cocycles
$Z^2_0(L^{n,m,p};L^{n,m,p})$ that vanish on the characteristic
vector $X_0$. Then, if $A=Z^2(L;L)\cap \mathrm{Hom}(L_0 \wedge
L_0,L_0)$ we have that

$$dim \ A=\left\{
\begin{array}{ll}
\displaystyle \frac{n(3n-2)}{8}& \mbox{if } n \mbox{ is even} \\
\\
\displaystyle \frac{3n^2-4n+1}{8}+ \lfloor \frac{n+1}{4} \rfloor &
\mbox{if } n \mbox{ is odd}
\end{array}
\right.$$

\begin{proof} It is convenient to distinguish the following
four cases where the reasoning for each case is not hard:
\begin{itemize}
\item[(1).] $n\equiv 0$ (mod 4).
 \item[(2).] $n\equiv 1$ (mod 4).
 \item[(3).] $n\equiv 2$ (mod 4).
  \item[(4).] $n\equiv 3$ (mod 4).
\end{itemize}
\end{proof}

\section{Dimension of $B=Z^2(L;L)\cap \mathrm{Hom}(L_0 \wedge
L_1,L_1) $}

\

In general, any cocycle $b \in Z^2(L;L)\cap \mathrm{Hom}(L_0
\wedge L_1,L_1) $ will be any skew-symmetric bilinear map from
$L_0 \wedge L_1$ to $L_1$ such that:

\

\begin{itemize}
\item[$(2)$]
$[X_i,b(X_j,Y_k)]-[X_j,b(X_i,Y_k)]-b([X_i,X_j],Y_k)+ $\\
$b([X_i,Y_k],X_j)+b(X_i,[X_j,Y_k])=0$ \hfill{$\forall \ X_i,X_j
\in L_0, \ Y_k \in L_1$}
\end{itemize}

\

\noindent with $b(X_0,X)=0$ $\forall X \in L$. This condition
reduces to

\begin{equation} \label{cocycle2}
[X_0,b(X_j,Y_k)]-b([X_0,X_j],Y_k)-b(X_j,[X_0,Y_k])=0,\ 1 \leq j
\leq n, \quad 1 \leq k \leq m
\end{equation}

\

In order to obtain the dimension of the space of cocycles $B$ we
apply an adaptation of the $\mathfrak{sl}(2,\CC)$-module Method
that we have already used in the precedent section.

\

Recall that if  $W_0,W_1,\dots,W_k$ are $\mathfrak{
sl}(2,\CC)$-modules, then the space $\mathrm{Hom}(\otimes_{i=1}^k
W_i,W_0)$ will be  a $\mathfrak{sl}(2,\CC)$-module in the
following natural manner:
$$
  (\xi \cdot \varphi)(x_1,\dots,x_k)=\xi \cdot \varphi(x_1,\dots,x_k)-
   \sum_{i=1}^{k}
    \varphi (x_1,\dots,\xi \cdot x_i,x_{i+1},\dots,x_n)
$$
with $\xi \in \mathfrak{ sl}(2,\CC)$ and $\varphi
\in\mathrm{Hom}(\otimes_{i=1}^k W_i,W_0)$. In particular, if $k=2$
and $V_0=W_1$, $V_1=W_2=W_0$, then
$$
    (\xi \cdot \varphi)(x_1,x_2)=\xi \cdot \varphi(x_1,x_2)-
 \varphi (\xi \cdot x_1,x_2)- \varphi (x_1,\xi \cdot x_2).
$$

An element $\varphi \in\mathrm{Hom}(V_0\otimes V_1,V_1)$ is said
to be invariant if
%$$
   $X_{+}\cdot \varphi =0$,
%$$
that is
\begin{equation} \label{maximal2}
      X_{+} \cdot \varphi(x_1,x_2)-
        \varphi (X_{+} \cdot x_1,x_2)- \varphi (x_1,X_{+} \cdot x_2)=0,
         \quad \forall x_1 \in V_0, \ \forall x_2 \in V_1.
\end{equation}

Note that $\varphi \in\mathrm{Hom}( V_0 \otimes V_1,V_1)$ is
invariant if and only if $\varphi$ is a maximal vector.

\

%On the other hand, we are going to consider the Lie superalgebra
%$L^{n,m}$ with basis $\{X_0,X_1,\dots,X_n,Y_1,\dots,Y_m\}$. By
%definition, a cocycle $\varphi$ belongs to $B_1$ it will be a
%skew-symmetric bilinear map:
%$$\varphi: L^{n,m}_0 \wedge L^{n,m}_1 \longrightarrow L^{n,m}_1$$
%such that $d \varphi=0$. That is, taking into account the law of
%$L^{n,m}$:
%\begin{equation}\label{cocycle} [X_0,\varphi(X_i,Y_j)]- \varphi
%([X_0,X_i],Y_j)-\varphi(X_i,[X_0,Y_j])=0, \mbox{ with } 1 \leq i
%\leq n, 1 \leq j \leq m \end{equation}

In this case we are going to consider the structure of irreducible
$\mathfrak{ sl}(2,\CC)$-module in
$V_0=<X_1,\dots,X_n>=L_0/\mathbb{C}X_0$ and in
$V_1=<Y_1,\dots,Y_n>=L_1$, thus in particular:
 $$\left\{ \begin{array}{ll} X_{+}\cdot X_i=X_{i+1}, & 1 \leq i \leq n-1\\[1mm]
X_{+}\cdot X_n=0&\\ [1mm] X_{+}\cdot Y_j=Y_{j+1}, & 1 \leq j \leq m-1\\[1mm]
X_{+}\cdot Y_m=0&
 \end{array}\right.$$

We identify the multiplication of $X_{+}$ and $X_i$ in the
$\mathfrak{ sl}(2,\CC)$-module  $V_0=<X_1,\dots,X_n>$, with the
bracket product $[X_0,X_i]$ in $L_0$. Analogously with $X_{+}
\cdot Y_j$ and $[X_0,Y_j]$. Thanks to these identifications, the
expressions (\ref{cocycle2}) and (\ref{maximal2}) are equivalent,
so we have the following result:

\begin{prop}Any skew-symmetric bilinear map $\varphi$,
$\varphi:V_0 \wedge V_1 \longrightarrow V_1$ will be an element of
$B$ if and only if $\varphi$ is a maximal vector of the
$\mathfrak{ sl}(2,\CC)$-module {\rm Hom} $(V_0 \wedge V_1,V_1)$,
with $V_0=<X_1, \dots,X_n>$ and $V_1=L_1$.
\end{prop}

%\begin{proof} If $\varphi \in B_1$, then $Im \varphi \in
%V_0$. In fact, if $X_0 \in Im \varphi$ it can be supposed that
%there exist $i,j$ such that $\varphi(Y_i,Y_j)=X_0$, but it
%constitutes a contradiction with the cocycle equation
%$$[x,\varphi(y,z)]-\varphi(z,[x,y])-\varphi(y,[x,z])=0,
%\quad \forall x\in L^{n,m}_0; y,z \in L^{n,m}_1$$ for $x=X_1$,
%$y=Y_i$ and $z=Y_j$. This proves the proposition and the equation
%\ref{cocycle}.
%\end{proof}

\begin{cor} As each $\mathfrak{ sl}(2,\CC)$-module has (up to nonzero scalar
multiples) a unique maximal vector, then the dimension of $B$ is
equal to the number of summands of any decomposition of {\rm Hom}
$(V_0 \wedge V_1,V_1)$ into direct sum of irreducible $\mathfrak{
sl}(2,\CC)$-modules.
\end{cor}

\

As each irreducible module contains either a unique (up to scalar
multiples) vector of weight $0$ or a unique vector of weight $1$,
then we have the following corollary.

%Thanks to the symmetric structure of the weights, instead of to
%sum the maximal vectors it is possible, and easier, to sum the
%vectors of weight 0 or 1.

\begin{cor}\label{cor2} The dimension of $B$ is equal to the dimension of the
subspace of  {\rm Hom} $(V_0 \wedge V_1,V_1)$ spanned by the
vectors of weight 0 or 1.
\end{cor}

\

Next,  we consider a natural basis of $\mathrm{Hom}(V_0 \wedge
V_1,V_1)$
%$$
%   {\mathcal B}=\{\varphi_{i,j}^s\}
%$$
%with
consisting of the following maps where $1 \leq s, j, l \leq m$ and
$1 \leq i, k \leq n$:
$$
   \varphi_{i,j}^s(X_k,Y_l)=\left\{
         \begin{array}{ll}
           Y_s & \mbox{ if }(i,j)=(k,l)\\
             0 & \mbox{ in all other cases}
         \end{array}
\right.
$$

\

 Thanks to
Corollary \ref{cor2} it will be enough to find the basis vectors
$\varphi_{i,j}^s$ with weight $0$ or $1$. It is not difficult to
see that the weight of an element $\varphi_{i,j}^s$ (with respect
to $H$) is
$$
  \lambda(\varphi_{i,j}^s)
    =\lambda(Y_s)-\lambda(X_i)-\lambda(Y_j)=n+2(s-i-j)+1.
$$

Thus, if $n$ is even then $\lambda(\varphi)$ is odd, and if $n$ is
odd then $\lambda(\varphi)$ is even. So, if $n$ is even it will be
sufficient to find the elements $\varphi_{i,j}^s$ with weight $1$
and if $n$ is odd it will be sufficient to find those with weight
$0$. To do that we consider the three sequences that correspond
with the weights of $V_0=<X_1,\dots,X_n>$,
$V_1=<Y_1,Y_2,\dots,Y_m>$ and $V_1=<Y_1,Y_2,\dots,Y_m>$:
$$
  -n+1,-n+3,\dots,n-3,n-1;
$$
$$
  -m+1,-m+3,\dots,m-3,m-1;
$$
$$
  -m+1,-m+3,\dots,m-3,m-1.
$$

We shall have to count the number of all possibilities to obtain
$1$ (if $n$ is even) or $0$ (if $n$ is odd). Remember that
$\lambda(\varphi_{i,j}^s)=\lambda(Y_s)-\lambda(X_i)-\lambda(Y_j)$,
where $\lambda(Y_s)$ belongs to the last sequence, and
$\lambda(X_i)$, $\lambda(Y_j)$ belong to the first and second
sequences respectively. Thus, we obtain the following theorem.

\

\noindent {\bf Theorem 2.} Let $Z^2(L;L)$ be the $2$-cocycles
$Z^2_0(L^{n,m,p};L^{n,m,p})$ that vanish on the characteristic
vector $X_0$. Then, if $B=Z^2(L;L)\cap \mathrm{Hom}(L_0 \wedge
L_1,L_1)$ we have that

$$dim \ B=\left\{
\begin{array}{ll}
\displaystyle \frac{4nm-n^2+1}{4}& \mbox{if } n \mbox{ is odd}, \ n < 2m+1 \\
\\
\displaystyle \frac{4nm-n^2}{4} & \mbox{if } n \mbox{ is even}, \
n < 2m+1 \\
\\
m^2 & \mbox{if } n \geq 2m+1
\end{array}
\right.$$

\begin{proof} It is convenient to distinguish the following
four cases where the reasoning for each case is not hard:
\begin{itemize}
\item[(1).] $n\equiv 0$ (mod 4).
 \item[(2).]$n\equiv 1$ (mod 4).
 \item[(3).]$n\equiv 2$ (mod 4).
  \item[(4).]$n\equiv 3$ (mod 4).
\end{itemize}
\end{proof}

\section{Dimension of $C=Z^2(L;L)\cap \mathrm{Hom}(L_0 \wedge
L_2,L_2) $}

\

Similarly to the previous section we can obtain the equivalent
result for $C$.

\

\noindent {\bf Theorem 3.} Let $Z^2(L;L)$ be the $2$-cocycles
$Z^2_0(L^{n,m,p};L^{n,m,p})$ that vanish on the characteristic
vector $X_0$. Then, if $C=Z^2(L;L)\cap \mathrm{Hom}(L_0 \wedge
L_2,L_2)$ we have that

$$dim \ C=\left\{
\begin{array}{ll}
\displaystyle \frac{4np-n^2+1}{4}& \mbox{if } n \mbox{ is odd}, \ n < 2p+1 \\
\\
\displaystyle \frac{4np-n^2}{4} & \mbox{if } n \mbox{ is even}, \
n < 2p+1 \\
\\
p^2 & \mbox{if } n \geq 2p+1
\end{array}
\right.$$

\section{conclusions}

The Theorems $1$, $2$ and $3$ together to those obtained in
\cite{dimension_color} and \cite{erratum2}, leads to obtain the
total dimension of the infinitesimal deformations of the model
$\ZZ_3$-filiform Lie algebra $L^{n,m,p}$. Thus, we have the
following theorem.

\

\noindent {\bf Main Theorem.} The dimension of the space of
infinitesimal deformations of the model $\ZZ_3$-filiform Lie algebra
$L^{n,m,p}$ that vanish on the characteristic vector $X_0$, is
exactly  $A+B+C+D+E+F$ where

\

\noindent$A=\left\{
\begin{array}{ll}
\displaystyle \frac{n(3n-2)}{8}& \mbox{if } n \mbox{ is even} \\
\\
\displaystyle \frac{3n^2-4n+1}{8}+ \lfloor \frac{n+1}{4} \rfloor &
\mbox{if } n \mbox{ is odd}
\end{array}
\right.$

\

\

\noindent$B=\left\{
\begin{array}{ll}
\displaystyle \frac{4nm-n^2+1}{4}& \mbox{if } n \mbox{ is odd}, \ n < 2m+1 \\
\\
\displaystyle \frac{4nm-n^2}{4} & \mbox{if } n \mbox{ is even}, \
n < 2m+1 \\
\\
m^2 & \mbox{if } n \geq 2m+1
\end{array}
\right.$

\

\

\noindent$C=\left\{
\begin{array}{ll}
\displaystyle \frac{4np-n^2+1}{4}& \mbox{if } n \mbox{ is odd}, \ n < 2p+1 \\
\\
\displaystyle \frac{4np-n^2}{4} & \mbox{if } n \mbox{ is even}, \
n < 2p+1 \\
\\
p^2 & \mbox{if } n \geq 2p+1
\end{array}
\right.$

\

$$D=\left\{
\begin{array}{llll}
\displaystyle \frac{m(m-1)}{2}& \mbox{\rm if } p \geq 2m-1& & \\
\\
\displaystyle \frac{1}{8}(4mp-p^2-2p-1) & \mbox{\rm if } p <
2m-1,&
p\equiv 1 \mbox{\rm (mod 4) and } m \ \mbox{\rm odd, or}\\
& & p\equiv 3 \mbox{\rm (mod 4) and } m \ \mbox{\rm even}\\ \\
\displaystyle \frac{1}{8}(4mp-p^2-2p+3) & \mbox{\rm if } p <
2m-1,&
p\equiv 3 \mbox{\rm (mod 4) and } m \ \mbox{\rm odd, or}\\
& & p\equiv 1 \mbox{\rm (mod 4) and } m \ \mbox{\rm even}\\ \\
\displaystyle \frac{1}{8}(4mp-p^2-2p) & \mbox{\rm if } p < 2m-1&
\mbox{\rm and } p \ \mbox{\rm even}
\end{array}
\right.$$

\

$$F=\left\{
\begin{array}{llll}
\displaystyle \frac{p(p-1)}{2}& \mbox{\rm if } m \geq 2p-1& & \\
\\
\displaystyle \frac{1}{8}(4pm-m^2-2m-1) & \mbox{\rm if } m <
2p-1,&
m\equiv 1 \mbox{\rm (mod 4) and } p \ \mbox{\rm odd, or}\\
& & m\equiv 3 \mbox{\rm (mod 4) and } p \ \mbox{\rm even}\\ \\
\displaystyle \frac{1}{8}(4pm-m^2-2m+3) & \mbox{\rm if } m <
2p-1,&
m\equiv 3 \mbox{\rm (mod 4) and } p \ \mbox{\rm odd, or}\\
& & m\equiv 1 \mbox{\rm (mod 4) and } p \ \mbox{\rm even}\\ \\
\displaystyle \frac{1}{8}(4pm-m^2-2m) & \mbox{\rm if } m < 2p-1&
\mbox{\rm and } m \ \mbox{\rm even}
\end{array}
\right.
$$

\

 $(1).$ If $m+p-n$ \mbox{ is even, then}

$$E =\left\{
\begin{array}{ll}
 mn & \mbox{if } p \geq m+n \\ \\
np-1 & \mbox{if } p < m+n, p=m-n+2 \\ \\
np & \mbox{if } p < m+n, p<m-n+2 \\ \\
\displaystyle \frac{1}{4}(-m^2-n^2-p^2+2np+2mn+2mp) & \mbox{if } p
< m+n, p>m-n+2, \\ & p\geq n-m+2 \\ \\
mp & \mbox{if } p < m+n, p>m-n+2, \\ & p< n-m+2 \\ \\
\end{array}
\right.$$

$(2).$ If $m+p-n$ \mbox{ is odd, then}

$$E =\left\{
\begin{array}{ll}
 mn & \mbox{if } p \geq m+n-1 \\ \\
np & \mbox{if } p < m+n-1, p \leq m-n+1 \\ \\
\displaystyle \frac{1}{4}(-m^2-n^2-p^2+2np+2mn+2mp+1) & \mbox{if }
p < m+n-1, p>m-n+1, \\ & p\geq n-m+1 \\ \\
mp & \mbox{if } p < m+n-1, p>m-n+1,\\ & p< n-m+1 \\ \\
\end{array}
\right.$$

\

\bibliographystyle{amsplain}

\end{document}